\documentclass{amsart}
\usepackage{amsmath,amssymb,amsthm,amsfonts,epsfig,topcapt}
\input{xy}                                                                      
\xyoption{all}
\CompileMatrices 

\newcommand{\com}[1]{
}
\newcommand{\Label}[1]{\label{#1}\com{#1}  }
\newcommand{\us}[1]{{\upshape{#1}}}
\newcommand{\bp}{\begin{pmatrix}}
\newcommand{\ep}{\end{pmatrix}}
\newcommand{\be}{\begin{equation}}
\newcommand{\ee}{\end{equation}}
\newcommand{\bs}{\begin{split}}
\newcommand{\es}{\end{split}}
\newcommand{\bc}{\begin{center}}
\newcommand{\ec}{\end{center}}
\newcommand{\ed}{{\rm d}}
\newcommand{\w}{{\mathchoice{\,{\scriptstyle\wedge}\,}{{\scriptstyle\wedge}}    
      {{\scriptscriptstyle\wedge}}{{\scriptscriptstyle\wedge}}}}                
\newcommand{\lhk}{\mathbin{\hbox{\vrule height1.4pt width4pt depth-1pt          
             \vrule height4pt width0.4pt depth-1pt}}}  

\newcommand{\liealgebra}[1]{{\mathfrak {#1}}}

\newcommand{\so}{\liealgebra{so}}

\newcommand{\liegroup}[1]{{\operatorname{#1}}}
\newcommand{\GL}{\liegroup{GL}}
\newcommand{\SO}{\liegroup{SO}}
\newcommand{\SU}{\liegroup{SU}}
\newcommand{\Un}{\liegroup{U}}
\newcommand{\Or}{\liegroup{O}}
\newcommand{\PSO}{\liegroup{PSO}}
\newcommand{\Spin}{\liegroup{Spin}}

\renewcommand{\Re}{\operatorname{Re}}
\renewcommand{\Im}{\operatorname{Im}}
\newcommand{\om}{\omega}

\newcommand{\gr}{\tilde{G}(2,7)}

\newcommand{\s}{S^6}

\newcommand{\R}{\mathbb R}
\newcommand{\C}{\mathbb C}
\newcommand{\Q}{\mathbb H}
\newcommand{\Z}{\mathbb Z}
\newcommand{\Oc}{\mathbb O}
\newcommand{\I}{\mathcal I}
\newcommand{\mcl}{\mathcal L}
\newcommand{\mcs}{\mathcal S}
\newcommand{\mcc}{\mathcal C}
\newcommand{\mcj}{\mathcal J}
\newcommand{\mck}{\mathcal K}

\newcommand{\ip}[2]{\langle {#1},\, {#2} \rangle}
\newcommand{\irr}{{\mathcal I}_{R}}

\newcommand{\p}{\varphi}
\newcommand{\G}{G_2}
\newcommand{\g}{\liealgebra{g}_2}

\newcommand{\csff}{\rm{I \negthinspace I}_{co}}

\newcommand{\e}{e_1 \w e_2}

\newcommand{\icr}{{\mathcal I}_{CR}}

\newcommand{\cp}[1]{{\C \mathbb{P}^{#1}}}

\newcommand{\db}{\bar{\partial}}

\newcommand{\La}{\Lambda}

\newcommand{\al}{\alpha}

\newcommand{\iso}{\to}
\newcommand{\ev}{{\bf e}}
\newcommand{\bx}{{\bf x}}
\newcommand{\ns}{\negthinspace}
\newcommand{\tr}{\; {}^t \ns}
\newcommand{\bi}{{\bf i }}
\newcommand{\bj}{{\bf j}}
\newcommand{\bk}{{\bf k}}
\usepackage[mathscr]{eucal}
\usepackage{bbm}
\usepackage{oldgerm}

\theoremstyle{plain}
\newtheorem{thm}{Theorem}[section]
\newtheorem{lem}[thm]{Lemma}

\newtheorem{prop}[thm]{Proposition}

\newtheorem{rem}[thm]{Remark}
\theoremstyle{definition}
\newtheorem{defn}{Definition}[section]

\newtheorem{exam}[thm]{Example}

\begin{document}
\title[2-Ruled Coassociative Cones]{Coassociative Cones That Are Ruled by 2-planes}
\author{Daniel Fox}
\address{Department of Mathematics\\ University of California Irvine\\ Irvine, California 92697}
\email{dfox@math.uci.edu}
\thanks{I would like to thank Robert Bryant for the help and encouragement he offered as my thesis adviser.  I would also like to thank Duke University, MSRI and Columbia University for their hospitality while I carried out my thesis work, upon which this article is based. This research was supported in part by NSF VIGRE Grant DMS-9983320 to Duke University.}
\date{July 1, 2005}
\begin{abstract}
It is shown that coassociative cones in $\R^7$ that are r-oriented and ruled by $2$-planes are equivalent to CR-holomorphic curves in the oriented Grassmanian of $2$-planes in $\R^7$.  The geometry of these CR-holomorphic curves is studied and related to holomorphic curves in $\s$.  This leads to an equivalence between associative cones on one side and the coassociative cones whose second fundamental form has an $\Or(2)$ symmetry on the other.  It also provides a number of methods for explicitly constructing coassociative $4$-folds.  One method leads to a family of coassociative $4$-folds whose members are neither cones nor are ruled by $2$-planes.  This family directly generalizes the original family of examples provided by Harvey and Lawson when they introduced coassociative geometry.
\end{abstract}
\maketitle

\tableofcontents

\begin{section}{Introduction}
In 1982 Harvey and Lawson introduced coassociative 4-folds in $\R^7$ as an example of a calibrated geometry \cite{hl}.  Using the perspective of exterior differential systems they proved a local existence theorem.  Harvey and Lawson offered a single finite dimensional family of explicit examples (see section \ref{sec:z3}).  In 1985 Mashimo classified the coassociative cones whose links are homogeneous \cite{mashimo}.  Then for almost twenty years the literature was quiet on the subject. The silence was broken in 2004 when techniques used for studying special Lagrangian geometry (which is another type of calibrated geometry) began to be applied to coassociative geometry \cite{lotayca}, \cite{ikm}, \cite{km}.  This renewal of activity was motivated by Joyce's construction of compact manifolds with $\G$ holonomy \cite{joyceg2} and the appearance of coassociative geometry in M-theory \cite{acharya}.  

In \cite{lotayca} Lotay studied coassociative 4-folds in $\R^7$ that are ruled by 2-planes.  He gave a local existence result and presented a method for using a holomorphic vector field on a naturally associated Riemann surface to deform a coassociative cone that is ruled by 2-planes.  His work is an extension of Joyce's work on ruled special Lagrangian 3-folds \cite{joyceruledsl}.

This article revisits the geometry of coassociative cones that are ruled by 2-planes using the methods introduced by Bryant in \cite{b1}.  This perspective realizes the Riemann surface that appeared in \cite{lotayca} as a CR-holomorphic for a $\G$-invariant CR-structure on $\gr$, the oriented Grassmanian of 2-planes in $\R^7$.  In fact, Proposition \ref{prop:cocr} asserts that r-oriented 2-ruled coassociative cones are equivalent to CR-holomorphic curves in $\gr$.  

Studying the geometry of CR-holomorphic curves leads to new methods for constructing coassociative cones that are ruled by 2-planes.  These methods arise from the close relationship between CR-holomorphic curves in $\gr$ and holomorphic curves\footnote{Some people may prefer using the term ``pseudo-holomorphic curves'' since the almost complex structure that is used on the ambient space is not integrable.} in $\s$.  This relationship can be interpreted as an equivalence between the associative cones and a certain family of coassociative cones.  The null-torsion holomorphic curves in $\s$ also form the backbone for a new family of coassociative 4-folds that directly generalizes the family introduced by Harvey and Lawson \cite{hl}.  The generic coassociative 4-fold in this new family is neither a cone nor does it admit a ruling by 2-planes, but it does retain an $S^1$-symmetry.   

Here is an outline of what follows.  Section \ref{sec:coassgeo} briefly introduces coassociative 4-folds.  Section \ref{sec:g2} provides an explicit description of $\g \subset \so(\R^7)$ that will be used for calculations. Section \ref{sec:eds} reviews the perspective of exterior differential systems (EDS). In Section \ref{sec:gstructures} the notion of a $G$-structure is reviewed.  The structure equations are written down for $\R^7 \times \G$ using the explicit description of $\g$ in Section \ref{sec:g2}.  Section \ref{sec:s6} reviews the $\SU(3)$-structure that is induced on $\s$ via the transitive $\G$-action and the resulting geometry of holomorphic curves.  This material will be central to describing the geometry of 2-ruled coassociative cones.  The treatment is based \cite{b2}.

Section \ref{sec:2ruled} describes the equivalence between r-oriented coassociative cones that are ruled by 2-planes and CR-holomorphic curves in $\gr$.  Section \ref{sec:crholcurv} addresses the geometry of CR-holomorphic curves in $\gr$.  In particular, two scalar invariants $a$ and $b$ and a holomorphic section $\rho$ of a holomorphic line bundle $\bar{\mcl}$ are invariantly associated to any CR-holomorphic curve in $\gr$.  The vanishing of each of the invariants is interpreted geometrically by relating CR-holomorphic curves in $\gr$ to holomorphic curves in $\s$ via the $\G$-equivariant map $p:\gr \to \s$.  When any of the three invariants vanishes there are methods for constructing explicit examples of such CR-holomorphic curves, and thus of coassociative cones.

Section \ref{sec:z3} describes the new family of coassociative 4-folds mentioned above.  Each member of it is constructed (using only differentiation) from a null-torsion holomorphic curve in $\s$.   
\end{section}
\begin{section}{Coassociative \us{4}-folds in $\R^7$}\Label{sec:coassgeo}
Coassociative 4-folds in $\R^7$ were introduced by Harvey and Lawson \cite{hl} as an example of a calibrated geometry.

\begin{defn}
Let $(X,g)$ be a Riemannian manifold.  A p-form $\psi$ on $X$ is said to be a \emph{calibration} if $\ed \psi =0$ and if for every point $x \in X$ and every orthonormal p-tuple of tangent vectors $v_1,\,\ldots,v_p$ in $T_xX$ the following inequality holds:
\be
\psi(v_1,\,\ldots,\,v_p)\leq 1.
\ee 
An orientable p-dimensional submanifold $f:M^p \to X$ is said to be \emph{calibrated} if 
\be
f^*(\psi)=\nu_f
\ee
where $\nu_f$ is the induced volume form on $M$.
\end{defn}
A calibrated manifold is a special type of minimal submanifold \cite{hl}.
\begin{thm}
A calibrated submanifold is volume minimizing in its homology class.
\end{thm}
The p-dimensional holomorphic submanifolds in a K\"ahler manifold $(X,\omega)$ are calibrated by $\frac{\omega^p}{p!}$.  It was K\"ahler geometry that motivated the introduction of calibrated geometry.

Harvey and Lawson introduced two new calibrated geometries in $\R^7$ \cite{hl}.  On $\R^7$ define the 3-form
\begin{equation}\Label{eq:3form}
\begin{split}
 \p=&\, \ed x_{5}\w \ed x_{6}\w \ed x_{7}
-\ed x_{5}\w(\ed x_{1}\w \ed x_{2}+\ed x_{3}\w \ed x_{4})\\
& -\ed x_{6}\w(\ed x_{1}\w \ed x_{3}+\ed x_{4}\w \ed x_{2})
-\ed x_{7}\w(\ed x_{1}\w \ed x_{4}+\ed x_{2}\w \ed x_{3}).
\end{split}
\end{equation}
If $\R^7$ is identified with the imaginary octonians $\Im(\Oc)$ and $x \cdot y$ represents octonianic multiplication of $x,y \in \Im(\Oc)$, then 
\be
\ip{x\cdot y}{z}=\p(x,y,z).
\ee
The affine linear stabilizer of $\p$ is $\R^7 \rtimes \G$, where $\G$ is the compact group associated to the 14-dimensional rank 2 exceptional Lie algebra $\g$.  In fact, $\G \subset \SO(7)$ and so the stabilizer of $\p$ also preserves the standard Euclidean inner product $\ip{\;}{}$ and the standard volume form $\nu=\ed x_1 \w \ldots \w \ed x_7$.  In \cite{hl} it is shown that both $\p$ and its Hodge dual $*\p$ are calibrations.  A 3-fold that is calibrated by $\p$ is called an \emph{associative} 3-fold. A 4-fold that is calibrated by $*\p$ is called a \emph{coassociative} 4-fold.  

There is an equivalent first order condition that is often easier to work with than the condition of being calibrated \cite{hl}.
\begin{lem}\Label{lem:coideal}
An immersed \us{4}-fold $f:M \to \R^7$ is coassociative if and only if $f^*(\p)=0.$
\end{lem}  Using this characterization Harvey and Lawson prove that locally a coassociative 4-fold depends on two functions of three variables (in the sense of exterior differential systems), showing that they are quite abundant locally.  
\end{section}
\begin{section}{An Explicit Description of $\g$}\Label{sec:g2}
An explicit description of $\g \subset \so(\R^7)$ will be essential for the calculations below.  The description used here is particularly suited for the study of coassociative geometry.  Let $T \subset \R^7$ denote the subspace $(x_1,x_2,x_3,x_4,0,0,0)$ and let $V^+ \subset \R^7$ denote the subspace  $(0,0,0,0,x_5,x_6,x_7)$.  Let  $\so(T)$ and $\so(V^+)$ be the natural subalgebras of $\so(\R^7)$.  There is an $\so(T)$-invariant decomposition 
\be
\La^2(T) \cong \La^2_{+}(T) \oplus \La^2_{-}(T),
\ee
where $\al \in \La^2_{\pm}(T)$ satisfies 
\be
\al \w \al =\pm |\al|^2 \ed x_1 \w \ed x_2 \w \ed x_3 \w \ed x_4.
\ee 
It is well known that $\so(T) \cong \so(\La^2_{-}(T)) \oplus \so(\La^2_{+}(T))$.  One can easily check that neither $\so(T)$ nor $\so(V^+)$ is a subalgebra of $\g$.  However, there is a copy of $\so(4)$ inside of their direct sum that is a subalgebra of $\g$.  Let $\so(3)^{\pm} := \so(\La^2_{\pm}(T))$ and let $\sigma^{\pm} : \so(T) \to \so(3)^{\pm}$.  The form $\p$ defines an isomorphism 
\be\Label{eq:normalbundle}
V^+ \; \tilde{\to} \; \La^2_{+}(T),
\ee
by the formula
\be
v \to v \lhk \p,
\ee
so let $\so(4) \subset \so(\R^7)$ act as $\so(T)$ on $T$ and as $\so(3)^+$ on $V^+$.  It turns out that this copy of $\so(4) \subset \g$ is the stabilizer of the coassociatve plane $T$.  It is useful to think of $T$ as the tangent space to a coassociative 4-fold and $V^+$ as its normal bundle.  Then the isomorphism in equation \eqref{eq:normalbundle} states that the the normal bundle of a coassociative 4-fold is isomorphic to its bundle of self-dual 2-forms.

Identify $T$ with $\Q$, the quaternions, via
$$(x_1,\,x_2,\,x_3,\,x_4,\,0,\,0,\,0) \to x_1 \cdot 1 +x_2 \cdot \bi + x_3 \cdot \bj + x_4 \cdot \bk \in \Q.$$  Let $\beta \in V^+ \otimes T$ and assume that it satisfies the equation
\begin{align}\Label{eq:betarel}
\bi \cdot \beta_5 + \bj \cdot \beta_6 + \bk \cdot \beta_7=0.
\end{align}
Then an element of $\g$ is equivalent to a pair $(\theta,\,\beta)$, where $\theta \in \so(T)$ and $\beta$ satisfies equation \eqref{eq:betarel}, via the map
\begin{align*}
(\theta,\,\beta) \to  \begin{pmatrix}\theta& -\tr{\beta}\\
\beta & \sigma^+(\theta) \end{pmatrix}.
\end{align*}  
To prove this one just needs to check that this defines a 14-dimensional Lie algebra that preserves $\p$.  This could be done by choosing a basis and computing.
\end{section}
\begin{section}{Exterior Differential Systems}\Label{sec:eds}
Lemma \ref{lem:coideal} indicates that exterior differential systems (EDS) arise naturally in coassociative geometry and in this article they will frequently be employed.  Exterior differential systems are not the height of fashion these days so I will briefly review the basic definitions.  For more information on EDS see the standard text \cite{bcggg}.  For a gentler introduction see \cite{il}.  

Let $X^n$ be a smooth manifold and let $\Omega(X)$ be the space of smooth sections of the bundle of differential forms, $\Omega(X)=\Omega^0(X)\oplus \Omega^1(X)\oplus \dots \oplus \Omega^n(X)$.  A differential ideal $\I$ on $X$ is an ideal contained in $\Omega(X)$ which is homogeneous and closed under differentiation.  Being homogeneous means that if $\alpha \in \I$, then $\alpha \cap \Omega^p(X) \in \I$ for all $p$.  Being closed under differentiation means that if $\alpha \in \I$ then $\ed \alpha \in \I$.  Together $(X,\I)$ will be referred to as an \emph{exterior differential system}, or EDS.

If one thinks of an EDS as an equation, the solutions are the \emph{integral manifolds}.  An immersed submanifold $f:\Sigma \to X$ is an integral manifold if $f^*(\I)=0$, i.e., $f^*(\alpha)=0$ for all $\alpha \in \I$.   If there is any ambiguity about what ideal is being refereed to I will say \emph{$\I$-integral manifolds}.  Lemma \ref{lem:coideal} implies that coassociative manifolds are the integral manifolds for the ideal generated by $\p$.  In general I will use the notation $\I=\langle \alpha,\,\dots,\omega \rangle$ to denote that $\I$ is generated by the finite set of forms $\alpha,\ldots, \omega$.

It is extremely useful to have good generators for an EDS $(X,\I)$.  Often such generators do not exist on the space $X$.  However, natural generators often do exist on a fiber bundle that has $X$ as its base.  This leads to the notion of a \emph{$G$-structure}.
\end{section}
\begin{section}{$G$-Structures}\Label{sec:gstructures}
A \emph{$G$-structure} on a smooth manifold $X$ is a principal subbundle $P$ of the coframe bundle $F$ that has structure group $G$.  The most common way of defining a $G$-structure on a manifold $X$ is as the subbundle of coframes in which the expressions for certain tensor fields on X (for example, the generators of an ideal $\I$) take a certain form.  Such a $G$-structure is known as a bundle of {\it adapted coframes}.  

Let $X$ be a smooth $n$-dimensional manifold.   A {\it coframe} of X at $x \in X$ is $ u_x \in T^*_xX \otimes \R^n$  s.t. $ u_x:T_x{X} \iso \R^n$ is an isomorphism.  Let $F$ denote the right principal $GL(n,\R)$-bundle of coframes and let $\pi: F \to X$ so that $\pi(u_x)=x$.  The right $GL(n,\R)$-action is defined by $u_x \cdot A:=A^{-1} \circ u_x$  for $A\in GL(n,\R)$.    

The frame bundle $F$ carries the tautological $\R^n$-valued 1-form $\om \in \Gamma(T^*F \otimes \R^n)$, defined by $\om_u := u \circ \pi_* $.  The tautological 1-form has the reproducing property: If $u$ is a local coframing of $X$, i.e., a local section of $\pi:F \to X$, then $u^*(\om)=u$.  There is also a natural bundle map $\ev:P \to TX \otimes (\R^n)^*$ defined by the equation $\pi_*=\ev \, \om$, or, in components, $\pi_* = e_i \om_i$.

Let $G \subset GL(n,\R)$ be a closed subgroup.  A $G$-structure $P$ on $X$ is a principal $G$-subbundle of $F$ .  The $\R^n$-valued 1-form $\om$ and $(\R^n)^*$-valued function $\ev$ pull back to $P$ and there $\om$ retains its reproducing property.

As an example take $X=\R^7$.  Then the coframe bundle is $F=\R^7 \times \GL(7,\R)$.  Let $U \subset \R^7$ be an open set and define a coframe $u:U \to F$ to be adapted if on $U$ 
\begin{align*}
\p &= u^{5}\w u^{6}\w u^{7}-u^{5}\w(u^{1}\w u^{2}+u^{3}\w u^{4})\\
&-u^{6}\w(u^{1}\w u^{3}-u^{2}\w u^{4})-u^{7}\w(u^{1}\w u^{4}+u^{2}\w u^{3}),
\end{align*}
when $\p$ is defined as in equation \eqref{eq:3form}.  Let $P$ be the subbundle of coframes that are adapted in this way.  Because $\G$ is the subgroup of $\GL(7,\R)$ that fixes $\p$, the principal bundle is $P\cong \R^7 \times \G$. When working on $\R^7$ (or $\R^n$) it is more common to let $\pi=\bx$ and to express $\pi_*=e_i\,\omega_i$ as $\ed \bx = e_i\,\omega_i$.  Then the equation $\ed e_i = e_j \omega_{ji}$ defines the skew-symmetric matrix of 1-forms $\omega_{ij}$ which, in this example, is just the Maurer-Cartan form on $\G$.  Taking $\ed$ of $\ed \bx = e_i \, \omega_i$ gives the structure equations involving $\ed \omega_i$.  Taking $\ed$ of the resulting equation then gives an expression for $\ed \omega_{ij}$.  Using the summation convention, these structure equations on $P=\R^7 \times \G$ are 
\be\Label{eq:streqn}
\begin{split}
\ed \bx &= e_i \omega_i\\
\ed e_i &= e_j\omega_{ji}\\
\ed \omega_i&= -\omega_{ij}\w\omega_j\\
\ed \omega_{ij}&=-\omega_{ik}\w\omega_{kj}
\end{split}
\ee
where $\omega_{ij}$ also satisfies $\omega_{ij}=-\omega_{ji}$ and
\be\Label{eq:g2rel}
\begin{split}
\omega_{67}&= \omega_{12}+\omega_{34}\\
\omega_{75}&= \omega_{13}+\omega_{42}\\
\omega_{56}&= \omega_{14}+\omega_{23}\\
\omega_{51}&=-\omega_{64}+\omega_{73}\\
\omega_{52}&=-\omega_{63}-\omega_{74}\\
\omega_{53}&= \omega_{62}-\omega_{71}\\
\omega_{54}&= \omega_{61}+\omega_{72}.
\end{split}
\ee
The algebraic identities satisfied by $\omega_{ij}$ are a consequence of the explicit description of $\g$ given in Section \ref{sec:g2}.
\end{section}
\begin{section}{Holomorphic Curves in $\s$}\Label{sec:s6}
Each coassociative cone that is ruled by 2-planes has a holomorphic curve at its core.  This holomorphic curve comes with a natural immersion into $\s$.  This will be described in Section \ref{sec:geoofinv}.  So that it is ready for use there, I will now review the geometry of holomorphic curves in $\s$ when the $\G$-invariant (non-integrable) complex structure is used.  This treatment is based on that given in \cite{b2}.

There is a unique irreducible action of $\G$ on $\R^7$.  The orbits consist of the origin and the 6-spheres centered at the origin.  Let  $\s$ denote the unit sphere in $\R^7$.  Choosing a point $s \in \s$ defines a map $u:\G \to \s$ by the rule $u(g)=g \cdot s$.  The stabilizer of a point is $\SU(3)$ and this makes $\G$ into a principal $\SU(3)$-bundle over $\s$. 

\hspace{1.5in} \xymatrix{\SU(3) \ar[r]& \G \ar[d]^{u}\\ & {\s}}\\
The point $s \in \s$ can be chosen so that $u=e_5$, where $e_5$ is a component of the adapted frame ${\bf e}$.  

The $\SU(3)$ stabilizer of a point in $\s$ acts on the tangent space $T_s \s \cong \R^6$ via the unique irreducible 6-dimensional representation.  Therefore $\G$, viewed as an $\SU(3)$-bundle over $\s$, is naturally a subbundle of the coframe bundle and thus the $\G$-action defines an $\SU(3)$-structure on $\s$.  To uncover the $\SU(3)$ structure equations, pullback the bundle $P=\R^7 \times \G$ to $\s$ using the inclusion $\iota:\s \to \R^7$ to get

\hspace{1.5in} \xymatrix{\G \ar[r]& \iota^{-1}(P) \ar[d]^{\bx}\\ & {\s}}\\
Then restrict to the $\SU(3)$-subbundle $P_{\SU(3)}$ for which $\bx=e_5$.  

\hspace{1.5in} \xymatrix{\SU(3) \ar[r]& P_{\SU(3)} \ar[d]^{e_5}\\ & {\s}}\\
Differentiating $\bx=e_5$ uncovers the identities $\omega_5=0$ and $\omega_i=\omega_{i5}$.  Thus the $\omega_{ij}$ provide a coframe of $P_{\SU(3)}$.  These identities allow the structure equations to be reorganized in a form that is suited to studying the geometry of holomorphic curves in $\s$.  Define the (complex valued) vector fields and 1-forms 
\begin{align*}
u&=e_5\\
(f_1,\,f_2,\,f_3)&=\frac{1}{2}(e_7+ie_6,\,-e_1-ie_2,\,-e_4+ie_3)\\
\begin{pmatrix} \theta_1 \\ \theta_2 \\ \theta_3  \end{pmatrix} &=\frac{1}{2} \begin{pmatrix} \omega_{65}+ i \omega_{75} \\ -\omega_{25}- i \omega_{15} \\ \omega_{35}- i \omega_{45} \end{pmatrix}\\
\alpha &=\begin{pmatrix} 
0&(\omega_{17}+\frac{1}{2}\omega_{35})&(-\omega_{36}-\frac{1}{2}\omega_{25})\\
-(\omega_{17}+\frac{1}{2}\omega_{35})&0&(-\omega_{23}+\frac{1}{2}\omega_{56}) \\
(\omega_{36}+\frac{1}{2}\omega_{25})&(\omega_{23}-\frac{1}{2}\omega_{56})&0
\end{pmatrix}\\
\beta &= \begin{pmatrix}
-\omega_{67}&(-\omega_{16}+\frac{1}{2}\omega_{45})&(-\omega_{37}+\frac{1}{2}\omega_{15})\\
(-\omega_{16}+\frac{1}{2}\omega_{45})&\omega_{12}&(-\omega_{13}-\frac{1}{2}\omega_{57}) \\
(-\omega_{37}+\frac{1}{2}\omega_{15})&(-\omega_{13}-\frac{1}{2}\omega_{57})&\omega_{34}
\end{pmatrix}\\
\kappa &=\alpha+i \beta 
\end{align*}
Then $f= (f_1,\;f_2,\;f_3)$ is an $\SU(3)$-adapted framing of $\s$.  The structure equations (equation \eqref{eq:streqn}) and $\g$ symmetries (equation \eqref{eq:g2rel}) now imply that
\begin{equation}\label{eq:hcstreqn}
\begin{split}
&\ed u =f (-2i\theta)+\bar{f}(2i\bar{\theta}) \\
&\ed f = u(-i\tr{\bar{\theta}}) + f \kappa - \bar{f} [\theta]\\
&\ed \theta = - \kappa \w \theta - [\bar{\theta}] \w \bar{\theta}\\
&\ed \kappa = -\kappa \w \kappa +3 \theta \w {}^t \bar{\theta} - {}^t
\theta \w \bar{\theta} \cdot Id\\
& \kappa =- \tr\,\bar{\kappa}, \;\;\;\; tr(\kappa)=0
\end{split}
\end{equation}
where $[a]:=\bp 0&a_3&-a_2\\-a_3&0&a_1\\a_2&-a_1&0 \ep$ for any vector $\tr{a}=(a_1,\;a_2,\;a_3)$.\footnote{Notice that these structure equations differ from those in \cite{b2} by two minus signs.  The map $(u,f,\theta,\kappa) \to (u,-f,-\theta,\kappa)$ gives the equivalence. }
  
The complex structure on $\s$ can be characterized as follows.  A 1-form $\alpha \in \Gamma(T^*\s \otimes \C)$ is declared to be a $(1,0)$-form if $u^*(\alpha) \in T^*\G \otimes \C$ is a $C^{\infty}(\G)$-linear combination of the $\theta_i$.  The $(1,1)$-form of the $\SU(3)$-structure is $\Omega=\frac{i}{2}(\theta_{1}\w\bar{\theta}_1+\theta_{2}\w\bar{\theta}_2+\theta_{3}\w\bar{\theta}_3)$ and the holomorphic volume form is $\Psi=\theta_1 \w \theta_2 \w \theta_3$. 

The complex structure allows one to define holomorphic curves.  Let $\phi:\Sigma^2 \to \s$ be an immersion of a (real) surface and let $\tilde{\phi}:\Sigma \to \G$ be a lift obtained by choosing an $\SU(3)$-adapted framing. 
\begin{defn}
$\phi:\Sigma \to \s$ is holomorphic if and only if $\tilde{\phi}^*(\theta_{i}\w\theta_{j})=0$ for all $i$ and $j$. 
\end{defn}

These holomorphic curves will play an important role in Sections \ref{sec:crholcurv} and \ref{sec:z3}.  Suppose that $\Sigma \subset \s$ is a holomorphic curve. When the tangent bundle of $\s$ is pulled back to $\Sigma$ it decomposes into three complex line bundles: the holomorphic tangent bundle $T^{1,0}\Sigma$, the first-normal bundle $N_1(\Sigma)$, and the second-normal (or bi-normal) bundle $N_2(\Sigma)$.  One can always adapt frames so that $T^{1,0}\Sigma =\C \cdot f_1$, $N_1(\Sigma) =\C \cdot f_3$, and $N_{2}(\Sigma) =\C \cdot f_2$ (Watch the indices!).  Let $P_{T^2}(\Sigma) \subset P_{\SU(3)}(\Sigma)=\phi^{-1}(\G)$ be the principal $T^2$-subbundle that preserves the decomposition $\phi^{-1}(T\s) = T^{1,0}\Sigma \oplus N_1(\Sigma) \oplus N_2(\Sigma)$. The functions and 1-forms $u$,$f$,$\theta$, and $\kappa$ pull back to $P_{T^2}(\Sigma)$ and I will use the same notation for their pull-backs.

The $\SU(3)$-adapted lifts of holomorphic curves in $\s$ are integral surfaces of a larger differential ideal.  The fact that $f_1$ is a framing of the tangent space implies that on $P_{T^2}(\Sigma)$ the relations $\theta_2=\theta_3=0$ hold.  Differentiating these relations leads to
\begin{align*}
0=\ed \theta_3 &= - \kappa_{31} \w \theta_1\\
0=\ed \theta_2 &= - \kappa_{21} \w \theta_1,\\
\end{align*}  
which imply that $\kappa_{31}=H_1\,\theta_1$ and $\kappa_{21}=\tilde{H}_1\,\theta_1$ for some complex valued functions $H_1$ and $\tilde{H}_1$.  The fact that $f_3$ spans the first normal bundle on $P_{T^2}(\Sigma)$ implies that $\tilde{H}_1=0$, i.e., that $\kappa_{21}=0$.  Differentiating this last condition results in the equation 
\begin{align*}
0=\ed \kappa_{21}&=-\kappa_{23} \w \kappa_{31}.
\end{align*}  
As long as $\kappa_{31}\neq0$ the last equation implies that $\kappa_{23}=H_2\,\theta_1$.  In \cite{b2} Bryant showed that the $H_1$ is a holomorphic function on $\Sigma$.  Therefore $\kappa_{31}$ vanishes identically or else $\kappa_{23}=H_2 \theta_1$.  He also showed that $\kappa_{31}=0$ implies that $\phi(\Sigma)$ is a round $S^2 \subset \s$ sitting in an associative 3-plane.  From now on assume that $\kappa_{31}\neq0$.  In summary, on $P_{T^2}(\Sigma)$ the following identities hold:
\begin{align*}
\theta_{2}&=0\\
\theta_{3}&=0\\
\kappa_{21}&=0\\
\kappa_{31}&=H_1 \theta_1\\
\kappa_{23}&=H_2 \theta_1,
\end{align*}  
and the last two are obtained by differentiating the first three.  

Let $\mcj=\langle \theta_{2},\,\theta_{3},\,\kappa_{21}\rangle$.  The argument above shows that the adapted lift of any holomorphic curve in $\s$ will be an integral surface.  A standard argument shows that any $\mcj$-integral surface $\Sigma \subset \G$ on which $\theta_1 \w \bar{\theta}_1\neq0$ is the adapted lift (so that $f_2$ spans $N_2(\Sigma)$) of the holomorphic curve $u:\Sigma \to \s$. 

In \cite{b2} Bryant defines a \emph{null-torsion} holomorphic curve to be one for which $\kappa_{23}=0$.  He shows that all of these holomorphic curves are naturally algebraic curves in the five-quadric $Q_5 \subset \cp{6}$ and that locally they can be explicitly described in terms of one arbitrary holomorphic function.  This will be described in more detail in Section \ref{sec:geoofinv}.  For now notice that the null-torsion holomorphic curves are equivalent to integral surfaces of the ideal $\mck=\langle \theta_{2},\,\theta_{3},\,\kappa_{21},\, \kappa_{23}\rangle$.
\end{section}
\begin{section}{Coassociative Cones Ruled by 2-planes}\Label{sec:2ruled}
One way to simplify the coassociative equations is to assume that the 4-fold admits a ruling.  In the context of calibrated geometry this technique has been applied to special Lagrangian 3-folds by Bryant \cite{b1} and Joyce \cite{joyceruledsl}.  Lotay has used this technique to study coassociative, Cayley and special Lagrangian 4-folds \cite{lotayca}, as well as associative 3-folds \cite{lotayass}.  He showed that this reduces the problem to studying an equation on a Riemann surface.  Complex geometry enters in a more elaborate way than he makes use of.  The geometry discussed below is similar to that found in \cite{b1}. 

\begin{defn}
$M^4 \subset \R^7$ admits a \emph{smooth ruling by $2$-planes} if there exists a smooth surface $\Sigma^2$ and a smooth map $\pi:M\rightarrow \Sigma^2$ such that for all $\sigma \in \Sigma$, $E_{\sigma}:=\pi^{-1}(\sigma)$ is a $2$-plane in $\R^7$.  Such a triple $(M,\pi,\Sigma)$ is said to be a \emph{$2$-ruled} 4-fold.  If there exists a continuous choice of orientation for $E_{\sigma}$ then $M$ is said to be \emph{r-oriented}.
\end{defn}

If $(M,\pi,\Sigma)$ is an r-oriented 2-ruled \emph{cone}, then all of the 2-planes $E_{\sigma}$ must pass through the origin.  This means that $\pi^{-1}$ can be viewed as a map $\tilde{\gamma}:\Sigma \to \gr$.  The triple $(M,\pi,\Sigma)$ will be referred to as \emph{nondegenerate} if $\gamma:\Sigma \to \gr$ is an immersion.

This suggests that the natural space for studying nondegenerate r-oriented 2-ruled 4-dimensional cones in $\R^7$ is the Grassmanian of oriented 2-planes in $\R^7$, $\gr$.
\begin{lem}\Label{lem:conesandsurfaces}
A nondegenerate r-oriented \us{2}-ruled \us{4}-dimensional cone in $\R^7$ is equivalent to a surface in $\gr$.
\end{lem}
\begin{proof}
That the cones lead to surfaces in $\gr$ has already been shown.  Given an immersed surface $\gamma:\Sigma \to \gr$, choose a local oriented orthonormal frame $(v_1,v_2)$ of the 2-plane and define $$\Gamma (r_1,r_2,\sigma ) = r_1 v_1 ( \sigma ) +r_2 v_2 ( \sigma ).$$  This is evidently an r-oriented 2-ruled cone.
\end{proof}

Let $\gamma:\Sigma \to \gr$ and let $\Gamma \left(r_1,r_2,\sigma \right)=r_1 v_1 \left(\sigma\right) +r_2 v_2 \left(\sigma \right)$ be the corresponding cone with image $M^4 \subset \R^7$. Let $\tilde{\Gamma}:\R^2 \times \Sigma \to P$ be a lift of $\Gamma$ so that $\e$ defines the ruling, i.e., $\tilde{\Gamma}(r_1,\,r_2,\,\sigma)=r_1\,e_1(\sigma)+r_2\,e_2(\sigma)$.  This is always possible because $\G$ acts transitively on $\gr$ \cite{hl}.  The condition for $\Gamma$ to be a coassociative immersion is 
\be
\Gamma^{*} \left(\varphi \right) =0.
\ee
By expanding in powers of $r_1,r_2,\ed r_1,\ed r_2$ this condition translates into the vanishing of two 1-forms and six 2-forms on $\Sigma$. These forms are the pullbacks under $\tilde{\Gamma}$ of the $\R^7 \rtimes \G$ invariant forms 
\begin{equation}\Label{eq:ideal1}
\ip{\ed e_i}{e_i \cdot e_j},\ip{\ed e_i}{\ed e_j \cdot e_k},
\end{equation}
where $i,j\in\{1,2\}$. Notice that these forms are defined on $\G$ and don't involve the $\R^7$ factor. This suggests defining $\irr$ to be the differential ideal on $\G$ generated by the 1-forms $\ip{\ed e_i}{e_i \cdot e_j}$ and the 2-forms $\ip{\ed e_i}{\ed e_j \cdot e_k} $ for $i,j \in \{1,\,2\}$.  Suppose that $\Sigma \subset \G$ is an integral surface of $\irr$ that projects to $\gr$ to be a surface.  Then the $\Gamma$ construction above leads to a 2-ruled 4-fold in $\R^7$, and the above calculation shows that it will be coassociative.  Thus the r-oriented 2-ruled coassociative cones are equivalent to integral surfaces of $(\G,\, \irr)$.

The integral surfaces of $(\G,\irr)$ all project to $\gr$ to be CR-holomorphic curves for a $\G$-invariant CR-structure.\footnote{For the bare essentials of almost-CR-structures the reader may want to consult \cite{b1}.}  To see this, first notice that the map $q:\G \to \gr$ is naturally given by $q:=\e$ and that because the stabilizer in $\G$ of an oriented 2-plane is $\Un(2)$, $q$ gives $\G$ the structure of a $\Un(2)$-bundle over $\gr$.  The 1-forms used as generators for $\irr$ are $\omega_{51}$ and $\omega_{52}$ and they span a subbundle of $T^* \G$ that is invariant under the stabilizer $\Un(2)$.  Therefore this bundle is the pull back of a subbundle of $T^*\gr$ whose annihilator $Q \subset T\gr$ is a real 8-plane bundle.  Therefore any $\irr$-integral surface must project to $\gr$ to be tangent to $Q$.  The ideal $\irr$ is itself $\Un(2)$-invariant.  This means that the $\Un(2)$-action preserves each fiber of $Q$ and so defines a complex structure $J$ on the vector bundle $Q \subset T\gr$.  The 2-forms in $\irr$ imply that a generic integral 2-plane is a complex line for that complex structure.  In fact:
\begin{prop}\Label{prop:cocr}
There is a complex structure $J$ on $Q \subset T\gr$ with the following properties:
\begin{itemize}
\item{$(Q,J)$ is a real analytic, Levi-flat almost CR-structure on $\G$ that is invariant under the $\G$ action.}
\item{Every CR-holomorphic curve gives rise to an r-oriented \us{2}-ruled coassociative cone via the $\Gamma$-construction.}
\item{Conversely, the surface in $\gr$ that is defined by any r-oriented \us{2}-ruled coassociative cone is a CR-holomorphic curve.}
\end{itemize}
\end{prop}
\begin{proof}
I'll begin by defining the almost CR-structure more explicitly.  The definition will make it clear that it is real analytic and the Levi-flatness will follow from the structure equations that it inherits from the Maurer-Cartan form on $\G$.  After this is done I will turn to the relationship between CR-holomorphic curves and the coassociative cones.  As always it will be convenient to calculate on $\G$ instead of $\gr$.  

On $\G$ define the complex valued 1-forms
\begin{equation}
\begin{split}
\zeta_3 = \omega_{31} + i \omega_{41},\;\zeta_4 = \omega_{32}+  i \omega_{42},\\
\zeta_6 = \omega_{61}- i \omega_{71},\; \zeta_7 = \omega_{62}- i \omega_{72}.
\end{split}
\end{equation}
These forms are semibasic\footnote{A differential form $\alpha$ on the total space of a fiber bundle $\pi:X \to B$ is \emph{semibasic} if at each $x\in X$ it is the pullback via $\pi^*_{x}$ of some form at $B_{\pi(x)}$.  This is equivalent to the vanishing of the contraction of $\alpha$ with any vector tangent to the fiber.} for the map $q:\G \to \gr$, as is the subbundle defined by $\omega_{51},\omega_{52}$.  In fact the real and imaginary parts of the $\zeta_i$ along with $\omega_{51}$ and $\omega_{52}$ for a basis for the $q$-semibasic forms on $\G$.  A 1-form $\alpha \in T^*\gr \otimes \C$ is defined to be of type $(1,0)$ if $q^{*}(\alpha)$ is a $C^{\infty}(\G)$ linear combination of the $\zeta_i$'s. In this way they define a complex structure $J$ on $Q$.

The $\G$ structure equations (Equations \eqref{eq:streqn} and \eqref{eq:g2rel}) imply the following structure equations for $\om_{51},\om_{52},\zeta_i$:
\be\Label{eq:crstreqn1}
\begin{split}
\ed \begin{bmatrix}\zeta_3\\ \zeta_4 \\ \zeta_6 \\ \zeta_7 \end{bmatrix} &\equiv \begin{bmatrix}i\omega_{34} &-\omega_{12}&\Phi&0\\\omega_{12}&i\omega_{34}&0&\Phi\\-\bar{\Phi}&0&-i(\omega_{12}+\omega_{34})&-\omega_{12}\\0&-\bar{\Phi}&\omega_{12}&-i(\omega_{12}+\omega_{34}) \end{bmatrix} \w \begin{bmatrix}\zeta_3\\ \zeta_4 \\ \zeta_6 \\ \zeta_7  \end{bmatrix}\\d\omega_{51} &\equiv \Re(\zeta_{3} \w \zeta_{7})+2\Im(\zeta_{3} \w \zeta_{6})+\Re(\zeta_{6} \w \zeta_{4})\\d\omega_{52} &\equiv 2\Re(\zeta_{4} \w \zeta_{7})+\Im(\zeta_{4} \w \zeta_{6})-\Im(\zeta_{7} \w \zeta_{3})
\end{split}
\ee
where $\Phi:=\omega_{63}+i\omega_{73}$ and the congruences are taken modulo $\omega_{51},\omega_{52}$.

 The equations in \eqref{eq:crstreqn1} show that the almost CR-structure defined by $\{ \omega_{51}, \omega_{52} , \zeta_{i}   \}$ is Levi-flat.  The almost CR-structure has codimension 2 because the fibers of $Q \subset T\gr$ are codimension 2; it is rank 4 because the fibers of $Q$ are 4-dimensional complex vector spaces.  

A real surface $\Sigma^2 \subset \gr$ is a CR-holomorphic if it is tangent to $Q$ and if $T_{\sigma}\Sigma \subset Q$ is a complex line for the complex structure $J$. Define the ideal
\be
{\mathcal I}_{CR}:=\langle  \omega_{51}, \omega_{52} , \zeta_{i} \w \zeta_{j} , \bar{\zeta_{i}} \w \bar{\zeta_{j}} \rangle.
\ee
The integral 2-planes for $\icr$ project under $q_*$ to be complex lines in $Q$ and so the integral surfaces that are transverse to $q$ project under $q$ to be the CR-holomorphic curves in $\gr$.  The CR-holomorphic curves for a Levi-flat almost-CR-structure of rank 4 depend on 6 functions of 1 variable, and so are at least locally abundant.

Define the 2-forms $\Upsilon_i$ by the equations 
\be
\begin{split}\Label{eq:crandr}
&\Upsilon_{1}+i\,\Upsilon_{2}=\zeta_{6} \w \zeta_{3} \\
&\Upsilon_{3}+i\,\Upsilon_{4}=\zeta_{7} \w \zeta_{4}\\
&\Upsilon_{5}+i\,\Upsilon_{6}=\zeta_{7} \w \zeta_{3}+ \zeta_{6} \w \zeta_{4}.  
\end{split}
\ee

Using the $\G$-structure equations one can check that 
\be\Label{eq:sigmas}
\begin{split}
\Upsilon_{1}&\equiv \frac{1}{2} \p(e_{1},\ed e_{1}, \ed e_{1})\\
\Upsilon_{2}&\equiv -\frac{1}{2}\p(e_{2},\ed e_{1}, \ed e_{1})\\
\Upsilon_{3}&\equiv \frac{1}{2}\p(e_{1},\ed e_{2}, \ed e_{2})\\
\Upsilon_{4}&\equiv -\frac{1}{2}\p(e_{2},\ed e_{2}, \ed e_{2})\\
\Upsilon_{5}&\equiv\p(e_{1},\ed e_{1}, \ed e_{2})\\
\Upsilon_{6}&\equiv-\p(e_{2},\ed e_{1}, \ed e_{2})\\
\end{split}
\ee
where the equivalences are modulo $\om_{51}$ and $\omega_{52}$.

Therefore  
\be
\irr = \langle \Upsilon_1,\, \ldots ,\, \Upsilon_{6},\,\om_{51},\,\om_{52}\rangle.
\ee
Using equation \eqref{eq:crandr}, this new characterization of $\irr$ implies that $\irr \subset \icr$.  Therefore every CR-holomorphic curve gives rise to a 2-ruled coassociative cone via the $\Gamma$-construction.

Now suppose that $\Gamma: \R^2 \times \Sigma \to \R^7$ is an r-oriented 2-ruled coassociative cone with image $M$.  Choose a lift $\tilde{\Gamma}:\R^2 \times \Sigma \to P$ so that $\tilde{\Gamma}(r_1,r_2,\sigma)=r_1 \,e_1 + r_2\,e_2$ and $e_1 \w e_2 \w e_3 \w e_4 = TM$.  Then by rotating the frame in the $\e$-plane it can be insured that $\om_{a1}=0$ for $a=5,6,7$.  This implies that $\zeta_6=0$.  First suppose that $\zeta_7 \neq 0$.  The fact that $\irr$ must vanish then implies that $\zeta_4=A\,\zeta_7$ and $\zeta_3 = B\,\zeta_7$, which implies that the corresponding surface $\gamma:\Sigma \to \gr$ is a CR-holomorphic curve.  

Now consider the case in which $\zeta_6=\zeta_7=0$.  This along with the vanishing of the 1-forms in $\irr$ imply that $\om_{a1}=\om_{a2}=0$ for $a=5,6,7$.  This is a linear equation on the second fundamental form of the coassociative cone.  It implies that the second fundamental form has an $\Or(2)$ symmetry.  In \cite{thesis} I show that it is an involutive condition and that the corresponding surface $\gamma:\Sigma \to \gr$ is a CR-holomorphic curve.  As will be described in Section \ref{sec:geoofinv}, this family is equivalent to the holomorphic curves in $\s$.
\end{proof}

\begin{rem}
Let $(x_1,\ldots,x_8)$ be standard coordinates on $\R^8 \cong \Oc$ so that $\Im(\Oc)=\{x_8=0\}$.  On $\Oc$ define the $4$-form $\Phi=\p \w \ed x_8+ *\p$.  The stabilizer of this $4$-form is $\Spin(7) \subset \SO(8)$.  Harvey and Lawson \cite{hl} show that $\Phi$ is a calibration. The calibrated $4$-folds are called \emph{Cayley} $4$-folds.  The Cayley cones that are ruled by $2$-planes are equivalent to holomorphic curves in $\tilde{G}(2,\Oc)$ for a non-integrable, $\Spin(7)$-invariant complex structure.  Proposition $\ref{prop:cocr}$ and its proof are actually special cases of the corresponding result for Cayley geometry.
\end{rem}

Proposition \ref{prop:cocr} along with the standard existence result for CR-holomorphic curves demonstrates that locally there are plenty of $2$-ruled coassociative cones, in fact six functions of one variable's worth of them.  This existence result is similar to Theorem 4.2 in \cite{lotayca}).

Proposition \ref{prop:cocr} shows that the 2-ruled condition reduces 4-dimensional coassociative geometry to a surface geometry.  It is this surface geometry that I will turn to now. 
\end{section}
\begin{section}{CR-Holomorphic Curves in $\gr$}\Label{sec:crholcurv}
\begin{subsection}{The Invariants $a$, $b$, and $\rho$}\Label{sec:crholcurvinv}
There are no $1^{st}$-order invariants of a coassociative 4-fold in $\R^7$ because $\G$ acts transitively on the Grassmanian of coassociative planes \cite{hl}.  There are many second order invariants.  A second order invariant for a 2-ruled coassociative 4-fold $(M,\pi,\Sigma)$ translates into a $1^{st}$-order invariant for the geometry of the corresponding CR-holomorphic curve $\Sigma \subset \gr$, and so one should expect $1^{st}$-order invariants for CR-holomorphic curves in $\gr$.  A more direct reason is that the structure group $\Un(2)$ does not act transitively on the complex lines in $Q_{v_1\w v_2} \cong \C^4$ for $v_1 \w v_2  \in \gr$.  

In this section I will introduce the scalar invariants $a,\,b \in \R$ for a CR-holomorphic curve in $\gr$.  There is a tautological complex line bundle $L$ defined on $\gr$.  It pulls back to be a holomorphic line bundle $\mathcal{L}$ over any CR-holomorphic curve and the geometry of the curve naturally defines a holomorphic section $\rho$ of its dual $\bar{\mcl}$.  Whether or not $\rho$ vanishes is another invariant of the CR-holomorphic curve.  The vanishing of each of these three invariants will be interpreted geometrically in Section \ref{sec:geoofinv}. 

To study the geometry of CR-holomorphic curves it will be useful to use the coframing introduced for holomorphic curves in $\s$. This coframing is related to the one introduced in the proof of Proposition \ref{prop:cocr} as follows:
\begin{equation}\Label{eq:s6cr2}
\begin{split}
2\theta_{1}&=- i \zeta_{3}+ \zeta_{4}\\
2\theta_{2}&=\omega_{52}+i\omega_{51}\\
2\theta_{3}&=i\zeta_{6} - \zeta_{7}\\
2\kappa_{21}&=i\zeta_6 +\zeta_7\\
2\kappa_{23}&=i\zeta_3+\zeta_4\\
\kappa_{11}&=-i\omega_{67}\\
\kappa_{22}&=i\omega_{12}\\
\kappa_{33}&=i\omega_{34}\\
\kappa_{31}&=-\bar{\Phi}=-(\omega_{63} - i \omega_{73}).
\end{split}
\end{equation}
The first equation in \eqref{eq:crstreqn1} can then be rewritten as
\be
\ed \begin{bmatrix} \theta_1 \\ \theta_3 \\ \kappa_{21} \\ \kappa_{23} \end{bmatrix} \equiv - \begin{bmatrix} 
\kappa_{11} & \kappa_{13}&0&0\\ 
\kappa_{31} & \kappa_{33}&0&0\\
0&0&\kappa_{22}-\kappa_{11}&-\kappa_{31}\\
0&0&-\kappa_{13}& \kappa_{22}-\kappa_{33}
\end{bmatrix} \w \begin{bmatrix} \theta_1 \\ \theta_3 \\ \kappa_{21} \\ \kappa_{23} \end{bmatrix}
\ee 
where the congruence is taken modulo $\om_{51}$ and $\om_{52}$, or modulo $\theta_2$ and $\bar{\theta}_2$.  

If $\Sigma \subset \G$ is the lift of a CR-holomorphic curve then locally  
\be
\begin{bmatrix} \theta_1 \\ \theta_3 \\ \kappa_{21} \\ \kappa_{23} \end{bmatrix} = \begin{bmatrix}A_1 \\ A_2 \\ B_1 \\ B_2 \end{bmatrix} \ed z
\ee
where $z$ is a holomorphic coordinate on $\Sigma$ and the $A_i$'s and $B_i$'s are complex functions.  A calculation shows that they are actually holomorphic functions.  Under a change of coframe given by $U \in \Un(2)$, $A$ and $B$ transform as
\begin{align*}
&A \to U A\\
&B \to (\det{\bar{U}}) \,\bar{U}B.
\end{align*}
These transformation properties imply that  
\begin{align}
a &= \tr{\bar{A}}A\\
b &= \tr{\bar{B}}B
\end{align}
are invariants for a CR-holomorphic curve.  If either $a$ or $b$ vanishes on a neighborhood then it vanishes identically since the $A_i$'s and $B_i$'s are holomorphic functions.

There is a tautological complex line bundle on $\gr$ that pulls back to CR-holomorphic curves to be holomorphic.  To see this let $v_1$ and $v_2$ be orthonormal vectors and recall that $\p$ defines a (noncommutative, nonassociative) multiplication on $\R^7$ by the rule 
\be
v_1 \cdot v_2 = \# \p(v_1,\,v_2,\,\,).
\ee
This multiplication has the property that $|v_1\cdot v_2|=|v_1||v_2|$.  Therefore $v_3:=v_2 \cdot v_1$ is a unit vector in $R^7$ and the real six plane perpendicular to it inherits a complex structure $J_{v_3}$ defined by $J_{v_3}(u)=u\cdot v_3$.  The real 2-plane $v_1 \w v_2$ is perpendicular to $v_3$ and turns out to be a complex line with respect to the induced complex structure on $v_3^{\perp}$.  This defines the tautological complex line bundle
\be
\C \to L \to \gr
\ee
whose fiber at $v_1 \w v_2 \in \gr$ is $\C \cdot (v_1+iv_2)$.  As is usual, $L$ is a subbundle of a trivial bundle, $L \subset \gr \times \R^7$.

Let $\gamma:\Sigma \to \gr$ be a CR-holomorphic curve and define $\mcl := \gamma^{-1}(L)$.  
\begin{lem}
The line bundle $\C \to \mcl \to \Sigma$ is holomorphic.
\end{lem}
\begin{proof}
A lift $\tilde{\gamma}:\Sigma \to \G$ can always be chosen so that $f_2$ is a unitary frame for $\mcl$.  This implies that $\sigma := \omega_1 - i\omega_2 $ along with a coframe of $\tilde{\gamma}(\Sigma)$ form a coframe of $\mcl$ and that $\omega_i=0$ except when $i=1,2$.  Then $\ed \sigma = -i \,\omega_{12} \w \sigma$.  This and the structure equations for $\gr$ make it clear that $\mcl$ has an integrable complex structure (since the $(1,0)$-forms algebraically generate their differential ideal) and that the map $\mcl \to \Sigma$ is holomorphic.  So $\mcl$ is a holomorphic line bundle over $\Sigma$. 
\end{proof}

Let $\bar{\mcl}$ be the dual of $\mcl$.  The transformation rules above imply that  
\be
\rho = \tr{B}A \bar{f}_2
\ee
is a well defined holomorphic section of $\bar{\mcl}$.

The vanishing of the invariants $a$, $b$, and $\rho$ can be interpreted according to the geometry of the image of $\gamma:\Sigma \to \gr$ under the $\G$-invariant map $p:\gr \to \s$.   
\end{subsection}
\begin{subsection}{Geometric Interpretations of the Invariants $a$, $b$, and $\rho$}\Label{sec:geoofinv}
There is a $\G$-equivariant map $p:\gr \to \s$ given by $p(v_1 \w v_2)=v_2 \cdot v_1$.  It fits into the commutative diagram of $\G$-equivariant maps

\hspace{2in} \xymatrix{ \G \ar[dd]^u \ar[dr]_q &\\& {\gr} \ar[dl]_p \\ {\s}&}\\

where $u=e_5$, $q=\e$, and $p \circ q=e_2 \cdot e_1 =e_5=u$.  The image in $\s$ of a CR-holomorphic curve is not arbitrary.
\begin{prop}\Label{prop:cohol}
Let $\tilde{\gamma}:\Sigma \to \G$ be the lift of a CR-holomorphic curve $\gamma:\Sigma \to \gr$ and let $\phi:=u \circ \tilde{\gamma}$.  Then $\phi:\Sigma \to \s$ is a holomorphic curve for the unique $\G$-invariant $\SU(3)$-structure. 
\end{prop} 
\begin{proof}
A surface in $\gr$ is CR-holomorphic if when it is lifted to $\G$ it is an integral surface for the ideal 
\be
\icr = \langle \om_{51},\,\om_{52},\,\zeta_i\w\zeta_j \rangle.
\ee
A surface in $\s$ is holomorphic if when it is lifted to $\G$ it is an integral surface for the ideal
\be\Label{eq:holcurvcond}
\mathcal{I}_{HC}=\langle \theta_{2}\w \theta_{3},\,\theta_{3}\w \theta_{1},\,\theta_{1}\w \theta_{2} \rangle.
\ee
The equations in \eqref{eq:s6cr2} show that $\I_{HC} \subset \icr$. 
\end{proof}
Proposition \ref{prop:cohol} begs the question, is every holomorphic curve in $\s$ in the image of a CR-holomorphic curve in $\gr$?  And, is there a way of producing a CR-holomorphic curve from a holomorphic curve in $\s$.  Both questions are about to be answered in the affirmative.  

There are three natural lifts of a holomorphic curve $\phi:\Sigma \to \s$ to $\gr$.  Each of the three complex line bundles $T\Sigma$, $N_1(\Sigma)$, $N_2(\Sigma)$ (see Section \ref{sec:s6}) has an oriented 2-plane in $\R^7$ as its fiber, and so each defines a lift $\Sigma \to \gr$.  It is the lift given by $N_2(\Sigma)$ that turns out to be useful here.  Following \cite{b2} this will be referred to as the \emph{binormal lift}.

For the rest of this section let $\gamma: \Sigma \to \gr$ be a CR-holomorphic curve and let $\tilde{\gamma}: \Sigma \to \G$ be a lift for which $\gamma=q \circ \tilde{\gamma}$.   

First consider the case in which $\rho=0$.
\begin{prop}
When $\rho=0$ but $a \neq 0$, $(\Sigma,\,\gamma)$ is the binormal lift of a holomorphic curve $\phi:\Sigma \to \s$.  
\end{prop}
\begin{proof}
Assume that $\rho=0$ but that $a \neq 0$.  This implies that $A$ and $B$ are orthogonal vectors in $\C^2$.  The $\Un(2)$-action can then be used to ensure that $\tr A=(A_1,\,0)$ and $\tr B = (0,\,B_2)$.  This implies that $\theta_3=\kappa_{21}=0$.  The fact that $\tilde{\gamma}(\Sigma)$ is an $\icr$-integral surface implies that $\theta_2=0$ and so $\tilde{\gamma}(\Sigma)$ is an integral surface of $\mcj=\langle \theta_2,\,\theta_3\,\kappa_{21}\rangle$.  This implies that it is an adapted lift of the holomorphic curve $\phi=u \circ \tilde{\gamma}:\Sigma \to \s$. In fact 
\be
\gamma=q \circ \tilde{\gamma} = \e=2if_2 \w \bar{f}_2
\ee
is the binormal lift of $(\Sigma,\phi)$.
\end{proof}

Now consider the special case in which $b=0$ but $a\neq0$.  

\begin{prop}\Label{prop:b=0}
When $b=0$ but $a\neq0$ $(\Sigma,\,\gamma)$ is the binormal lift of a \emph{null-torsion} holomorphic curve in $\s$.
\end{prop}
\begin{proof}
The situation is the same as in the last proof except that now the extra condition $b=0$ implies that $\kappa_{23}=0$.  This is the extra condition needed to make $(\Sigma,\,\tilde{\gamma})$ an integral surface of $\mck=\langle \theta_2,\,\theta_3,\,\kappa_{21},\,\kappa_{23}\rangle$, and thus the lift of a null-torsion holomorphic curve.
\end{proof}

In \cite{b2} Bryant showed that the null-torsion holomorphic curves are naturally algebraic curves in the five-quadric, $Q_5$.  It is well known that $Q_5 \cong \gr$ and the binormal lift of a holomorphic curve $\phi:\Sigma \to \s$ naturally lifts $\Sigma$ to $Q_5$.  He showed that the null torsion condition makes the lift of $\Sigma$ a holomorphic curve with respect to the standard holomorphic structure on $Q_5 \subset \cp{6}$.  Not only that, but it is also tangent to a non-integrable holomorphic 2-plane field.  Bryant showed how these holomorphic curves can locally be expressed in terms of a single arbitrary holomorphic function, and that globally every compact Riemann surface admits an immersion as such a holomorphic curve in $Q_5$.  By Proposition \ref{prop:b=0} this leads to a way of explicitly constructing coassociative cones. 

The coassociative cones constructed from null-torsion holomorphic curves in this way are naturally K\"ahler surfaces \cite{thesis}.  In fact, these cones can be deformed so that they are no longer conical but they remain 2-ruled and K\"ahler.  This should correspond to the deformation method introduced by Lotay in \cite{lotayca}.\footnote{At least locally there is an even larger family of deformations that preserve the property of being 2-ruled.  In \cite{thesis} I show that this larger family depends on a function $f$ satisfying $\Delta f = 2 f$ on $\Sigma$ and on a section of an affine complex line bundle that satisfies a $\db$-type equation.  Unfortunately this larger class of deformations destroys the K\"ahler structure.  The case in which $f=0$ should reduce to Lotay's construction using holomorphic vector fields.} 

Finally consider the case in which $a=0$ but $b\neq0$.
\begin{prop}
When $a=0$ but $b\neq0$ the image of $(\Sigma,\,\gamma)$ in $\s$ is a single point and $\Sigma$ is naturally an algebraic curve in $\cp{2}$.
\end{prop}
\begin{proof}
If $a=0$ then $\theta_i=0$ for $i=1,\,2,\,3$, which implies that $u \circ \tilde{\gamma}\;(\Sigma)=s_0$ for some point $s_0 \in \s$.  Therefore $\gamma(\Sigma)$ is contained in the fiber of $p:\gr \to \s$ which is diffeomorphic to $\cp{2}$.  On the fibers of $p$ the structure equations in \eqref{eq:hcstreqn} reduce to
\begin{align*}
\ed \kappa &= -\kappa \w \kappa\\
\tr\,\bar{\kappa}&=-\kappa\\
tr(\kappa)&=0
\end{align*}
which shows that the fiber inherits the standard symmetric space structure of $\cp{2}$.  The fact that $\tilde{\gamma}(\Sigma)$ is an integral of $\icr$ implies that $\kappa_{21} \w \kappa_{23} =0$, which shows that it will project down to $p^{-1}(s_0)\cong \cp{2}$ to be a holomorphic curve, and thus an algebraic curve. 
\end{proof}

The results so far can be summarized as
\begin{thm}\Label{thm:crgeoint}
The vanishing of the invariants have the following geometric interpretations:
\begin{enumerate}
\item{$\rho=0$ and $a \neq 0$ $\iff$ $\Sigma$ is the binormal lift of a holomorphic curve in $\s$.}
\item{$b=0$ and $a \neq 0$ $\iff$ $\Sigma$ is the binormal lift of a null-torsion holomorphic curve in $\s$ and is thus an algebraic curve in $Q_5$.}  
\item{$a=0$ and  $b \neq 0$ $\iff$ $u:\gr \to \s$ restricts to $\Sigma$ to be constant.  In this case $\Sigma \subset \gr$ is an algebraic curve in the fiber $u^{-1}(p_0) \cong \cp{2}$.}
\end{enumerate} 
\end{thm}

In \cite{thesis} I showed that the CR-holomorphic curves for which $\rho=0$ but $a\neq0$ correspond to the coassociative cones whose second fundamental forms have an $\Or(2)$ stabilizer.  Using Theorem \ref{thm:crgeoint} this gives an equivalence between the holomorphic curves in $\s$ and the coassociative cones whose second fundamental form has an $\Or(2)$ stabilizer.  It is well known that the cone on a holomorphic curve in $\s$ is an associative \us{3}-fold \cite{b2}.  This leads to the following correspondence:
\begin{thm}\Label{thm:asscoass}
There is a bijection between the associative cones and the coassociative cones whose second fundamental forms have an $\Or(2)$-stabilizer. 
\end{thm}

In \cite{ktw} Kong, Terng and Wang show that the equations describing holomorphic curves in $\s$ form an integrable system.  They provide a method for constructing explicit examples of holomorphic curves that have an $S^1$-symmetry.  Theorem \ref{thm:asscoass} shows that these methods will also provide explicit examples of 2-ruled coassociative cones.

The coassociative cones that are equivalent to the family of CR-holomorphic curves for which $a=0$ are also K\"ahler surfaces.  The fact that the image of $\Sigma$ in $\s$ is a point implies that the corresponding coassociative 4-fold is contained in $s_0^{\perp} \cong \R^6 \subset \R^7$.  It is well known that the $G_2$-structure on $\R^7$ induces an $\SU(3)$-structure on any hyperplane and that a coassociative 4-fold contained in a hyperplane is just a complex surface for that (flat) $\SU(3)$-structure.  If the coassociative 4-fold is ruled by 2-planes then part three of Theorem \ref{thm:crgeoint} implies that as a complex surface in $\C^3$ it is ruled by complex lines, and so is equivalent to a holomorphic curve in $\cp{2}$.  

Finally, I want to point out that Theorem \ref{thm:crgeoint} naturally breaks up the geometry of 2-ruled coassociative cones into two components.  The first component is a holomorphic curve $\Sigma \subset \s$.  Then the bundle $\cp{2} \to \gr \to \s$ pulls back to be a holomorphic $\cp{2}$-bundle over $\Sigma$.  The second component is a holomorphic section of this holomorphic $\cp{2}$-bundle. 

The geometric descriptions given here offer a number of ways that one could apply existing methods to construct explicit examples of r-oriented 2-ruled coassociative 4-folds.  The relationship with null-torsion holomorphic curves also provides a method for constructing explicit coassociative 4-folds that are not ruled by 2-planes.
\end{subsection}
\end{section}
\begin{section}{A Surface Bundle Construction}\Label{sec:z3}
When Harvey and Lawson introduced coassociative geometry \cite{hl} they offered a single explicit nontrivial family of coassociative 4-folds in $\R^7$.  
\begin{exam}{(Harvey and Lawson)}\Label{hlexample}
Let $(x_1,\,\ldots ,\, x_7)$ be the standard coordinates on $\R^7$. Let $\SU(2)$ act on $\R^7 \cong \R^4 \oplus \R^3$ in the standard irreducible ways on $\R^4\cong \C^2$ and $\R^3$ so that $(x_1,\,x_2,\,x_3,\,x_4)$ are the standard coordinates on $\R^4$ and $(x_5,\,x_6,\,x_7)$ are the standard coordinates on $\R^3$.  Let $\R^2 \subset \R^7$ be the $(x_1,x_5)$ subspace and let $\mcc$ by the algebraic curve defined by the equation $x_5(x_5^2-\frac{5}{4}x_1^2)^2=k^5$, where $k \in \R$.  When $\SU(2)$ acts on $\R^7$, the curve $\mcc$ sweeps out a 4-fold, $M$.  Harvey and Lawson show that $M$ is coassociative.
\end{exam}

The coassociative 4-folds constructed by Harvey and Lawson naturally have the structure of a surface bundle.  Harvey and Lawson point out that the $\SU(2)$ orbits that foliate this 4-fold are graphs of the Hopf fibration $S^1 \to S^3 \to S^2$, where the base space is a round 2-sphere contained in the (associative) 3-plane $0 \oplus \R^3 \subset \R^4 \oplus \R^3$.  In fact the whole 4-fold $M$ is a surface bundle $\mcs \to M \to S^2$.  Each fiber $\mcs$ is a surface of revolution whose geodesic circles are the fibers of the Hopf map and whose profile curve is the algebraic curve $\mcc$ from Example \ref{hlexample}.    

The surface bundle structure of $M$ suggests a way to generalize this family.  One could try to maintain the surface bundle structure while deforming the base space.  As the base space is deformed the fibers must be dragged along with it.  To do this one must relate the geometry of the fibers to that of the base.  This can be done by viewing the base as a holomorphic curve in $\s$.  As described in Section \ref{sec:s6}, each holomorphic curve $\phi:\Sigma \to \s$ has a tangent bundle, a normal bundle, and a binormal bundle.  The $S^3$ orbit described in Example \ref{hlexample} is the principal $\Un(1)$-bundle associated to the second normal bundle of the holomorphic curve $S^2 \subset \s$.  The fiber $\mcs$ of $M$ is contained in $\R \cdot \phi \oplus N_2(S^2)$.  Its axis of revolution is $\R \cdot \phi$.  When $S^2$ is replaced by any other \emph{null-torsion} holomorphic curve $\Sigma \subset \s$, the 4-fold consisting of the analogous surface bundle over $\Sigma$ is still coassociative.

\begin{exam}{(Surface Bundle)}\Label{exam:surfacebundles}
Let $\phi:\Sigma \to \s$ be a connected {\it null-torsion} holomorphic curve, let $P_{T^2}(\Sigma)$ be the bundle of adapted frames in which $f_1$ spans $T\Sigma$, $f_2$ spans $N_1(\Sigma)$ and $f_3$ spans $N_2(\Sigma)$.  (This is a different adaptation than was used in section \ref{sec:geoofinv}.)  Let $E=\R \cdot \phi \oplus N_{2}(\Sigma)$. This is a rank 3 real vector bundle over $\phi(\Sigma)$.  In each fiber $E_{\sigma}$ let $(w,\,z,\,\theta)$ be cylindrical coordinates s.t. $\frac{\partial}{\partial w}=\phi$, $z$ is the radial coordinate perpendicular to the axis of rotation defined by $\phi$, and $\theta$ is the angular coordinate.  Using these coordinates define the surface of rotation $\mcs_{\sigma} \subset E_{\sigma}$ by the equation   
\be
w(w^2-\frac{5}{4}z^2)^2=k^{5}.
\ee
The profile of this surface is shown in Figure \ref{profilecurve} along with with the asymptotes.  The surface $\mcs_{\sigma}$ has two connected components.  Let $\pi:M_{\Sigma} \to \Sigma$ be the surface bundle over $\phi(\Sigma) \subset \s$ whose fiber is $\mcs_{\sigma}$.
\begin{thm}\Label{thm:z3geom}
$M_{\Sigma}$ is a coassociative $4$-fold with two connected components corresponding to the two connected components of $\mcs_{\sigma}$. When $\Sigma$ is closed and $k\neq0$ $M_{\Sigma}$ is complete. When $k \to 0$ the fibers $\mcs_{\sigma}$ degenerate into the union of $N_2(\Sigma)_{\sigma}$ and a round cone obtained by rotating the line $w=\frac{\sqrt{5}}{2}z$ about the $\R \cdot \phi$ axis.  The corresponding surface bundles are coassociative cones.
\end{thm}
\end{exam}
\begin{figure}
\begin{center}
\topcaption[The (Algebraic) Profile Curve]{The (algebraic) profile curve of $\mcs_{\sigma}$.  The vertical axis is the axis of rotation so that $\mcs_{\sigma}$ has two connected components.  The asymptotes are the lines $w=\pm \frac{\sqrt{5}}{2}z$ and $w=0$.}
\framebox{
\epsfig{file=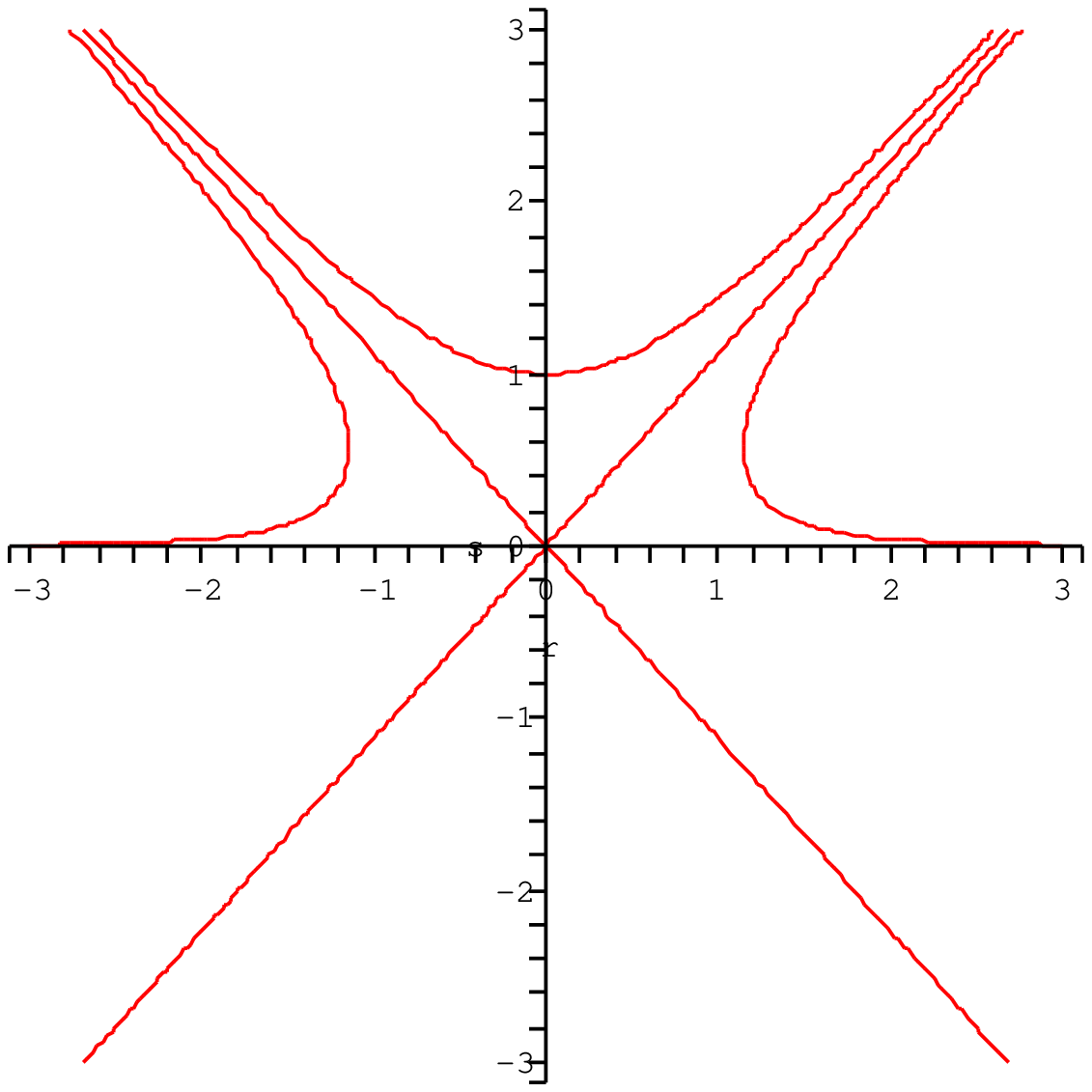, width=2 in}}
\Label{profilecurve}                                                           
\end{center}  
\end{figure}
\begin{proof}
It is enough to show that the tangent planes of $M_{\Sigma}$ are coassociative planes.  I will compute a framing $(h_1,\,h_2,\,h_3,\,h_4)$ of $M_{\Sigma}$ and check that $\p(h_i,\,h_j,\,h_k)=0$ for $1 \leq i,j,k \leq 4$.  

Let $P_{T^2}(\Sigma) \subset \G$ be as in Example \ref{exam:surfacebundles}.  Recall that $u:P_{T^2}(\Sigma) \to \s$ and that $u \circ \tilde{\phi}=\phi$ where $\tilde{\phi}:\Sigma \to P_{T^2}(\Sigma)$ is an adapted lift.  Define
\be
\tilde{\bx} = w \,u + 2\, z\,  \Im(f_3):P_{T^2}(\Sigma) \times \R^2 \to \R^7.
\ee
The image of $\bx$ is the 3-plane bundle $E$.  When $(z,w)$ are pulled back to the profile curve $\mcc$ the image of $\tilde{\bx}$ restricts to be $M_{\Sigma}$.  Let $t=\frac{w}{z}$ so that 
\begin{align*}
z&=k\,t^{-\frac{1}{5}}(t^2-\frac{5}{4})^{-\frac{2}{5}}\\
w&=k\,t^{\frac{4}{5}}(t^2-\frac{5}{4})^{-\frac{2}{5}}
\end{align*}
is a parametrization of $\mcc$.  Then restricting $(z,w)$ to the profile curve defines a map
\be
\bx: P_{T^2}(\Sigma) \times \R \to M
\ee
which is given explicitly by 
\be
\bx = w(t) \,u + 2\, z(t)\,  \Im(f_3)=k\,t^{-\frac{1}{5}}\,(t^2-\frac{5}{4})^{-\frac{2}{5}}\,[t\,u+2\Im(f_3)].
\ee
when $t \neq 0,\, \pm \frac{\sqrt{5}}{2}$.
This description suggests two natural legs for a framing of $M_{\Sigma}$.  One leg is given by the unit vector field $h_1 = |\frac{\ed \bx}{\ed t}|^{-1} \frac{\ed \bx}{\ed t}$, which is tangent to the profile curve.  The other is the vector field tangent to the $S^1$-symmetry corresponding to the fact that $\mcs_{\sigma}$ is a surface of revolution. This is $h_2 = -2 \Re(f_3)$.  So far the partial framing is
\begin{align*}
h_1&=\left( \frac{z^2-4w^2}{z^2+4w^2}\right) 2\Im(f_3)+ \left(\frac{-4zw}{z^2+4w^2} \right)u\\
h_2&=-2\Re(f_3).
\end{align*}

Two more legs are needed to have a framing of $M_{\Sigma}$.  Calculating $\ed \bx$ with the structure equations for a holomorphic curve (Equation \eqref{eq:hcstreqn}) gives
\begin{align*} 
\ed \bx  = &h_1\, \eta_1 + h_2\, \eta_2 \\
& +[4\,w\,\Re(f_1) -2z\,\Re(f_2)]\Im(\theta_1)+[4\,w\,\Im(f_1)+2z\, \Im(f_2)]\Re(\theta_1)
\end{align*}
for some 1-forms $\eta_1$ and $\eta_2$ whose precise formulas won't be needed.  This suggests defining
\begin{align*}
h_3&=|4\,w\,\Im(f_1)+2z\,\Im(f_2)|^{-\frac{1}{2}}[4\,w\,\Im(f_1)+2z\,\Im(f_2)]\\
h_4&=|4\,w\,\Re(f_1)-2z\,\Re(f_2)|^{-\frac{1}{2}}[4\,w\,\Re(f_1)-2z\,\Re(f_2)].
\end{align*}
This leads to the framing
\begin{align*}
h_1&=\left( \frac{z^2-4w^2}{z^2+4w^2}\right) 2\Im(f_3)+ \left(\frac{-4zw}{z^2+4w^2} \right)u\\
h_2&=-2\Re(f_3)\\
h_3&= \frac{4\,w\,\Im(f_1)+2z\,\Im(f_2)}{\sqrt{z^2+4w^2}}\\
h_4&= \frac{4\,w\,\Re(f_1)-2z\,\Re(f_2)}{\sqrt{z^2+4w^2}}.\\
\end{align*}
Then requiring that  
\begin{align*}
\ed \bx&=h_1\,\eta_{1}+h_2\,\eta_{2}+h_3\,\eta_{3}+h_4\,\eta_{4}\\
\end{align*}
defines $\eta_3$ and $\eta_4$, although their exact formulas will not be used either.  What is important is that this shows that $(h_1,h_2,h_3,h_4)$ is a framing of $TM_{\Sigma}$. 

Now recall that the framing $(u,f)$ was defined in terms of the $\G$-adapted framing $(e_1,\,\ldots,\,e_7)$ with dual coframing $(\omega_1,\ldots,\,\omega_7)$ in which\begin{align*}
\p &= \omega_{5}\w \omega_{6}\w \omega_{7}-\omega_{5}\w(\omega_{1}\w \omega_{2}+\omega_{3}\w \omega_{4})\\
& -\omega_{6}\w(\omega_{1}\w \omega_{3}+\omega_{4}\w \omega_{2})-\omega_{7}\w(\omega_{1}\w \omega_{4}+\omega_{2}\w \omega_{3}).
\end{align*}
By expressing the $h_i$'s in terms of the $e_i$'s it is easy to check that $\p(h_i,h_j,h_k)=0$ for $i,j,k \leq 4$.  Therefore $M_{\Sigma}$ is coassociative.
\end{proof}

This family first arose from studying coassociative 4-folds whose second fundamental form has a nontrivial stabilizer \cite{thesis}. The space of coassociative second fundamental forms is isomorphic to $\csff=\R_-^3 \otimes S_0^2(\R^3_+)$, which is an irreducible representation for $\PSO(4) \cong \SO(3)^- \times \SO(3)^+$.  Requiring that a coassociative second fundamental form have a nontrivial stabilizer forces it to lie in a proper subspace of $\csff$, which is a second order equation on a coassociative 4-fold. A coassociative 4-fold has the geometry described in Example \ref{exam:surfacebundles} if and only if the stabilizer of its second fundamental form contains the diagonal copy of $\Z_3 \subset \SO(3)^- \times \SO(3)^+$.  
\end{section}
\bibliographystyle{amsplain}
\bibliography{coassociativecones}
\end{document}